# New Ore-type Conditions for Hamilton Cycles and Spanning Trees with few leaves

Zhora Nikoghosyan

**Abstract**

New Ore-type conditions for Hamilton cycles and spanning trees with few leaves are presented.

## 1. Introduction

We consider only finite, undirected graphs with no loops or multiple edges. A good reference for any undefined terms is [1]. Let $G = (V(G), E(G))$ be a graph with vertex set $V(G)$ and edge set $E(G)$. For a vertex subset $X \subseteq V(G)$, we use $|X|$ to denote the cardinality of $X$. For a vertex $v \in V(G)$, we denote by $N(v)$ the neighborhood of $v$ in $G$, $d(v) = |N(v)|$ the degree of $v$ in $G$.

A path (simple path) of order $m$, denoted by $P_m$, is a sequence of different vertices $v_1, \ldots, v_m$, denoted by $v_1 v_2 \ldots v_m$, such that $v_{i-1} v_i \in E(G)$ for all $2 \leq i \leq m$. Similarly, a cycle (simple cycle) of order $m$ is a sequence of distinct vertices $v_1, \ldots, v_m$, denoted by $v_1 v_2 \ldots v_m v_1$, such that $v_{i-1} v_i \in E(G)$ and $v_m v_1 \in E(G)$ for all $2 \leq i \leq m$. In particular, for $m = 2$, $v_1 v_2 v_1$ is a cycle of order 2; and for $m = 1$, $v_1 v_1$ is a cycle of order 1. So, by the definition, every vertex (edge) can be interpreted as a cycle of order 1 (2, respectively). A graph $G$ is Hamiltonian, if it contains a Hamilton cycle - a simple spanning cycle. This extension of the cycle concept allows us to avoid the undesirable additional condition $n \geq 3$ in almost all Hamiltonian sufficient conditions for graphs on $n$ vertices.

A leaf of a tree is a vertex with degree one. A tree with at most $k$ leaves is called an $k$-ended tree, where $k \geq 2$ is an integer.

Suppose that all the edges of a path $P = v_1 v_2 \ldots v_m$ are directed from $v_1$ to $v_m$. For $x \in \{v_2, v_3, \ldots, v_m\}$, we denote the predecessor of $x$ on $P$ by $x^-$. For $U \subseteq V(P) \setminus \{x_1\}$, we denote $U^- = \{u^- : u \in U\}$.

In 1952, Dirac [3] gave the first sufficient condition for a graph to be Hamiltonian in terms of minimum degree $\delta$.

**Theorem A** [3]. Every graph on n vertices with $\delta \geq \frac{n}{2}$ is Hamiltonian.

In 1961, Ore [5] relaxed the minimum degree condition in Theorem A in terms of degree sums of nonadjacent vertices.

**Theorem B** [5]. Every graph on $n$ vertices is Hamiltonian, if $d(u) + d(v) \geq n$ for each nonadjacent vertices $u$ and $v$.

In 1976, an analogous condition was given for $k$-ended spanning trees [4].

**Theorem C** [4]. For a given integer $k \geq 2$, every connected graph on $n$ vertices satisfying $d(u) + d(v) \geq n - k + 1$ for each nonadjacent vertices $u, v$, has a $k$-ended spanning tree.

Theorem C can be considered as a starting point for all analogous results concerning spanning trees with few leaves and branch vertices (a branch vertex of a tree is a vertex with degree more than two). The PDF-file of paper [4] is available at http://ysu.am/files/1-1599118662-.pdf. Theorem C is also mentioned in [6]. However, the literature mainly refers to the work of Broersma and Tuinstra [2], published in 1998 (39 citations in Scopus).

In this paper we present two new Ore-type conditions for Hamilton cycles and $k$-ended spanning trees by relaxing the conditions in Theorems B and C on special induced subgraphs in terms of $K_m$ (complete graph), $K_{m,n}$ (complete bipartite graph), $P_m$ (simple path on m vertices) and $W + e$ (adding an edge $e$ to a subgraph $W$).

**Theorem 1**. Let G be a graph on $n \geq 4$ vertices satisfying $d(x) + d(y) \geq n$ for each nonadjacent vertices $x, y$ of every induced $K_{1,2} \cup K_1$, $K_3 \cup K_1$, $K_{1,3}$, $K_{1,3} + e$ and $P_4$. If $G$ has a vertex $v$ with $d(v) \geq 2$, then $G$ is Hamiltonian.

**Theorem 2**. For a given integer $k \geq 2$, every connected graph on $n$ vertices satisfying $d(x) + d(y) \geq n - k + 1$ for each nonadjacent vertices $x, y$ of every induced $K_{1,2} \cup K_1$, $K_3 \cup K_1$, $K_{1,3}$, $K_{1,3} + e$ and $P_4$, has a spanning $k$-ended tree.

Observe that the graphs $K_3 \cup K_1$, $K_{1,3}$ and $P_4$ can be alternatively obtained from $(K_{1,2} \cup K_1) + e$. Then the following two corollaries follow immediately.

**Corollary 1**. Let $G$ be a graph on $n \geq 4$ vertices satisfying $d(x) + d(y) \geq n$ for each nonadjacent vertices $x, y$ of every induced $K_{1,2} \cup K_1$, $(K_{1,2} \cup K_1) + e$ and $K_{1,3} + e$. If $G$ has a vertex $v$ with $d(v) \geq 2$, then $G$ is Hamiltonian.

**Corollary 2**. For a given integer $k \geq 2$, every connected graph on $n$ vertices satisfying $d(x) + d(y) \geq n - k + 1$ for each nonadjacent vertices $x, y$ of every induced $K_{1,2} \cup K_1$, $(K_{1,2} \cup K_1) + e$ and $K_{1,3} + e$, has a spanning $k$-ended tree.

## 2. Proofs

**Proof of Theorem 1**. Assume first that $G$ is not connected and let $H_1, H_2, \ldots, H_h$ be the connected components of $G$. By the hypothesis, there is a vertex $v$ such that $d(v) \geq 2$. Assume without loss of generality that $v \in V(H_1)$. Then we can form a path $Q = xvy$ in $H_1$. Let $z \in V(H_2)$ and let $F$ be a graph induced on $\{x, y, z, v\}$. If $xy \notin E(G)$, then clearly $F = K_{1,2} \cup K_1$. Otherwise, we have $F = K_3 \cup K_1$. In both cases, we have $vz \notin E(F)$. By the hypothesis, $d(v) + d(z) \geq n$, implying that

$$n \geq |V(H_1)| + |V(H_2)| \geq (d(v) + 1) + (d(z) + 1) \geq n + 2,$$

a contradiction. So, we can assume that $G$ is connected.

Let $P = x_1 x_2 \ldots x_p$ be a longest path in $G$.

**Case 1.** $x_1 x_p \in E(G)$.

Put $C = x_1 x_2 \ldots x_p x_1$. If $p = n$, then $C$ is a Hamilton cycle in $G$ and we are done. Otherwise, since $G$ is connected, there is an edge $y_1 y_2 \in E(G)$ such that $y_1 \in V(C)$ and $y_2 \notin V(C)$. Choose a vertex $y_3 \in V(C)$ such that $y_1 y_3 \in E(C)$. Then $C - y_1 y_3 + y_1 y_2$ is a path longer than $P$, a contradiction.

**Case 2.** $x_1 x_p \notin E(G)$.

**Case 2.1.** Either $d(x_1) \geq 2$ or $d(x_p) \geq 2$.

Assume without loss of generality that $d(x_1) \geq 2$. Recalling that $P$ is a longest path in $G$, we have $N(x_1) \subset V(P)$. Then we can choose $i \geq 3$ such that $x_1 x_i \in E(G)$.

**Case 2.1.1.** $x_{i-1} x_p \in E(G)$.

Put $C_1 = x_1 x_2 \ldots x_{i-1} x_p x_{p-1} \ldots x_i x_1$. If $p = n$, then $C_1$ is a Hamilton cycle in $G$. If $p \neq n$, then we can argue as in Case 1.

**Case 2.1.2.** $x_{i-1} x_p \notin E(G)$.

Since $i \geq 3$, we have $x_{i-1} \neq x_1$, implying that $x_1, x_{i-1}, x_i, x_p$ are distinct vertices. Let $R$ be a graph induced on $\{x_1, x_{i-1}, x_i, x_p\}$. If $x_1 x_{i-1} \notin E(G)$ and $x_i x_p \notin E(G)$, then $R = K_{1,2} \cup K_1$. Next, if $x_1 x_{i-1} \notin E(G)$ and $x_i x_p \in E(G)$, then $R = K_{1,3}$. Further, if $x_1 x_{i-1} \in E(G)$ and $x_i x_p \notin E(G)$, then $R = K_3 \cup K_1$. Finally, if $x_1 x_{i-1} \in E(G)$ and $x_i x_p \in E(G)$, then $R = K_{1,3} + x_i x_p$. Since $x_1 x_p \notin E(R)$, by the hypothesis, $d(x_1) + d(x_p) \geq n$. If $N^-(x_1) \cap N(x_p) \neq \emptyset$, then we can argue as in Case 2.1.1. Let $N^-(x_1) \cap N(x_p) = \emptyset$. Observing that $x_p \notin N^-(x_1) \cup N(x_p)$, we obtain

$$p \geq |N^-(x_1)| + |N(x_p)| + |\{x_p\}| \geq d(x_1) + d(x_p) + 1 \geq n + 1,$$

a contradiction.

**Case 2.2.** $d(x_1) = d(x_p) = 1$.

Let $R$ be a graph induced by $\{x_1, x_2, x_3, x_p\}$. If $p \geq 5$, then $R = K_{1,2} \cup K_1$. If $p = 4$, then $R = P_4$. Now let $p = 3$. Since $n \geq 4$ and $G$ is connected, we have $yx_2 \in E(G)$ for some vertex $y \notin P$. Since $P$ is a longest path in $G$, $\{y, x_1, x_3\}$ is an independent set of vertices. Hence $R = K_{1,3}$. By the hypothesis, $2 = d(x_1) + d(x_3) \geq n$, a contradiction. Theorem 1 is proved. ∎

**Proof of Theorem 2.** Let $P = x_1 x_2 \ldots x_p$ be a longest path in $G$.

**Case 1.** $p \geq n - k + 2$.

Since $G$ is connected, we can form a spanning tree $T$ in $G$ including $P$ as a subpath with end vertices $x_1$ and $x_p$. Clearly, $T$ has at most

$$|V(G)| - |V(P) \setminus \{x_1, x_p\}| = n - p + 2 \leq k$$

leaves and we are done.

**Case 2.** $p \leq n - k + 1$.

If $p = n$, then $P$ is a Hamilton path in $G$. This means that $G$ has an $r$-ended spanning tree with $r = 2 \leq k$ and we are done. Now let $p \neq n$. If $x_1 x_p \in E(G)$, then recalling that $G$ is connected, we can form a path longer than $P$ as in proof of Theorem 1 (Case 1), a contradiction. So, we can assume that $x_1 x_p \notin E(G)$.

**Case 2.1.** Either $d(x_1) \geq 2$ or $d(x_p) \geq 2$.

Assume without loss of generality that $d(x_1) \geq 2$. Recalling that $P$ is a longest path in $G$, we have $N(x_1) \subset V(P)$. Choose $i \geq 3$ such that $x_1 x_i \in E(G)$. If $x_{i-1} x_p \in E(G)$, then using the cycle $x_1 x_2 \ldots x_{i-1} x_p x_{p-1} \ldots x_i x_1$ and recalling that $p \neq n$, we can form a path longer than $P$ as in proof of Theorem 1 (Case 1). Let $x_{i-1} x_p \notin E(G)$.

Clearly, $x_1, x_{i-1}, x_i, x_p$ are distinct vertices. Let $R$ be a graph induced on $\{x_1, x_{i-1}, x_i, x_p\}$. As in proof of Theorem 1 (Case 2.1.2), we have either $R = K_{1,2} \cup K_1$ or $R = K_{1,3}$ or $R = K_3 \cup K_1$ or $R = K_{1,3} + x_i x_p$. By the hypothesis, $d(x_1) + d(x_p) \geq n$. If $N^-(x_1) \cap N(x_p) \neq \emptyset$, then we can argue as in Case 2.1.1. Let $N^-(x_1) \cap N(x_p) = \emptyset$. Observing that $x_p \notin N^-(x_1) \cup N(x_p)$, we obtain

$$p \geq |N^-(x_1)| + |N(x_p)| + |\{x_p\}| \geq d(x_1) + d(x_p) + 1 \geq n + 1,$$

a contradiction.

**Case 2.2.** $d(x_1) = d(x_p) = 1$.

**Case 2.2.1.** $2 \leq p \leq 3$.

If $p = 2$, then clearly $n = p = 2$ and $P$ is a 2-ended spanning tree in $G$. Let $p = 3$. If $n = 3$, then again P is a 2-ended spanning tree in $G$. Now let n$\geq$ 4. Since $G$ is connected and $d(x_1) = d(x_p) = 1$, we have $yx_2 \in E(G)$ for some vertex $y \notin V(P)$. Since $P$ is a longest path in $G$, $\{y, x_1, x_3\}$ is an independent set of vertices and $R = K_{1,3}$. By the hypothesis, $2 = d(x_1) + d(x_3) \geq n \geq 4$, a contradiction.

**Case 2.2.2.** $p \geq 4$.

Let $R$ be a graph induced by $\{x_1, x_2, x_3, x_p\}$. If $p = 4$, then $R = P = P_4$. If $p \geq 5$, then $R = K_{1,2} \cup K_1$. By the hypothesis, $2 = d(x_1) + d(x_p) \geq n \geq 4$, a final contradiction. Theorem 2 is proved. ∎

Institute for Informatics and Automation Problems of NAS RA
P. Sevak 1, Yerevan 0014, Armenia
E-mail: zhora@iiap.sci.am; ORCID ID: 0000-0001-9043-0458